\documentclass [a4paper,twoside,11pt]{article}

\usepackage[latin1]{inputenc}
\usepackage{amsmath,amsfonts,amssymb}
\usepackage{vmargin,graphicx,theorem}
\usepackage[english]{babel}
\usepackage{enumerate}

\setpapersize[portrait]{A4}
\setmarginsrb{2.5cm}{2.5cm}{2.5cm}{2cm}{0cm}{0.1cm}{0.5cm}{2cm}
\newcommand{\disp}{\displaystyle}
\newcommand{\dN}{\ensuremath{\mathbb{N}}}

\newcommand{\dR}{\ensuremath{\mathbb{R}}}

\newtheorem{ethm}{Theorem}[section]

\newtheorem{eprop}[ethm]{Proposition}

\newtheorem{elem}[ethm]{Lemma}

\newtheorem{edefi}[ethm]{Definition}

\newtheorem{erem}[ethm]{Remark}

\newcommand{\proofend}{~$\rhd$}
\newcommand{\proofbegin}{~$\lhd$}

\newenvironment{eproof}
               {\noindent {\emph{\textbf{Proof}}}\\\proofbegin~}
               {\proofend\\}

\newcommand{\p}[4]{{#3}\!\left#1{#4}\right#2}

\newcommand{\ABS}[1]{\ensuremath{{\left| #1 \right|}}} % |1|
\newcommand{\PAR}[1]{\ensuremath{{\left(#1\right)}}} % (1)
\newcommand{\NRM}[1]{\ensuremath{{\left\Vert #1\right\Vert}}} % ||1||
\renewcommand{\phi}{\varphi}
\renewcommand{\geq}{\geqslant}

\newcommand{\varf}[1]{\mathbf{Var}_{#1}}
\newcommand{\entf}[1]{\mathbf{Ent}_{#1}}

\newcommand{\ent}[2]{\p(){\entf{#1}}{#2}}

\newcommand{\var}[2]{\p(){\varf{#1}}{#2}}

\def\disp{\displaystyle}

\newcommand{\A}{\ensuremath{\mathcal A}}

\newcommand{\GI}{\mathbf{L}}
\newcommand{\PT}[1]{\mathbf{P}_{\!#1}}
\newcommand{\Pt}[1][t]{\ensuremath{\mathbf{P_{\!#1}}}}

\begin{document}

\title{\sl Logarithmic Sobolev inequality for diffusion semigroups }
\author{ Ivan Gentil\thanks{Ceremade, UMR CNRS 7534, Universit\'e Paris-Dauphine.}  }
\date{\today}
\maketitle\thispagestyle{empty}
%\abstract{Course given in Grenoble at the  2009 summer school :
 %Optimal transportation, Theory and applications}
\begin{abstract}
Through the main example of the Ornstein-Uhlenbeck semigroup, the Bakry-Emery criterion is presented as a main tool to get  functional inequalities as Poincar\'e or logarithmic Sobolev inequalities. Moreover an alternative method using the optimal mass transportation, is also given to obtain the logarithmic Sobolev inequality.
\end{abstract}


\smallskip


%\subjclass{Primary 35B40, 35K10, 60J60.}

\bigskip

%\noindent
%{\bf Key words:} 

\section{Introduction}
The goal of this course is to introduce inequalities as Poincar\'e or logarithmic Sobolev  for diffusion semigroups. We will focus more on examples than on the general theory.  A main tool to obtain those inequalities is the so-called Bakry-Emery $\Gamma_2$-criterium. This criterium is well known to prove such inequalities and has been also  used many times  for other problems, see for instance \cite{bakryemery,bakrytata}. We will focus on the example of the Ornstein-Uhlenbeck semigroup and on the $\Gamma_2$ -criterium. 

\medskip

In section~\ref{sec-ou} we investigate the  main example of the Ornstein-Uhlenbeck semigroup whereas in section~\ref{curvature} we show how the $\Gamma_2$-crierium implies such inequalities.  In section~\ref{sec-tr}, we will explain an alternative method to get a logarithmic Sobolev inequality under curvature assumption. It is called the {\it mass transportation method} and has been introduced recently, see \cite{c,ov,cnv,villani}.  By this way we will also obtain an another inequality called the   {\it Talagrand inequality or $\mathcal T_2$ inequality}.

\section{The Ornstein-Uhlenbeck semigroup and the Gaussian measure}
\label{sec-ou}
In the general setting if $(X_t)_{t\geq0}$ is a Markov process on $\dR^n$ then the family of operators : 
$$
\Pt (f)(x)=E(f(X_t)), 
$$
where $X_0=x$ and a smooth function $f$, defined is Markov semigroup on $\dR^n$.  There are two main examples. The first one is the heat semigroup which is associated to the Brownian motion on $\dR^n$. In this course we will study the second one which is the Ornstein-Uhlenbeck semigroup. We will see that the Ornstein-Uhlenbeck semigroup is associated to a linear stochastic differential equation driven by a Brownian motion. 

\medskip

In this note a smooth function $f$ in $\dR^n$ is a function such that  all computation done as integration by parts are justified, for example $\mathcal C_c^\infty(\dR^n)$. 

\subsection{Definition and general properties}
\begin{edefi}
Let define the family of operator $(\Pt)_{t\geq0}$ :  if  $f\in\mathcal C_b(\dR^n)$ then  
\begin{equation}
\label{eq-defou}
\Pt f(x)=\int f(e^{-t}x+\sqrt{1-e^{-2t}}y)d\gamma(y),
\end{equation}
where 
$$
d\gamma(y)=\frac{e^{-{\ABS{y}^2}/{2}}}{\PAR{2\pi n}^{n/2}}dy
$$ is the standard Gaussian distribution in $\dR^n$ and $\ABS{\cdot}$ is the Euclidean  norm on $\dR^n$. 

The family of operator $(\Pt)_{t\geq0}$ is called the Ornstein-Uhlenbeck semigroup. 
\end{edefi}
\begin{erem}
Let $(X_t)_{t\geq 0}$ be  a Markov process, solution of the stochastic differential equation 
\begin{equation}
\label{eq-eds}
\left\{
\begin{array}{l}
dX_t=\sqrt{2}dB_t-X_tdt\\
X_0=0.\\
\end{array}
\right. 
\end{equation}
Since the stochastic differential  equation is linear, there is an explicit solution
 $$
 X_t=e^{-t}X_0+\int_0^t \sqrt{2}e^{s-t}dB_s,
 $$ 
and equation~\eqref{eq-defou} is known as the Mehler Formula.  Moreover It\^o's formula gives that for all continuous and bounded functions $f$ on $\dR^n$ 
$$
\Pt f(x)=E_x(f(X_t)).
$$
 \end{erem}

\begin{eprop}
The Ornstein-Uhlenbeck semigroup is a linear operator satisfying the following properties : 
\begin{enumerate}[(i)]
\item $\PT{0}=Id$
\item For all functions  $f\in\mathcal C_b(\dR^n)$, the map $t\mapsto \Pt f $ is continuous from $\dR^+$ to $\mathcal L^2(d\gamma)$.
\item For all $s,t\geq 0$ one has $\Pt\circ\PT{s}=\PT{s+t}$.
\item $\Pt 1=1$ and  $\Pt f\geq 0$ if $f\geq 0$. 
\item $\NRM{\Pt f}_\infty\leq\NRM{f}_\infty$.
\end{enumerate}
We say that the Ornstein-Uhlenbeck semigroup is a {\it Markov semigroup} on $(\mathcal C_b(\dR^n),\NRM{\cdot}_\infty)$.
\end{eprop}
\begin{eproof}
We will give only some indications of the proof. First it is easy to prove items $(i)$, $(ii)$, $(iv)$ and $(v)$. 

For the item $(iii)$, you just have to compute the Ornstein-Uhlenbeck as follow : 
$\Pt f(x)=E(f(e^{-t}x+\sqrt{1-e^{-2t}}Y))$ where $Y$ is a random variable with a Gaussian distribution. Then compute $\Pt(\PT{s} f)$ to obtain $\PT{t+s} f$.  
In fact, since the solution of the stochastic differential equation~\eqref{eq-eds} is  a Markov process then $(iii)$ is a natural property of the Ornstein-Uhlenbeck semigroup. 
\end{eproof}

\begin{eprop}
\label{prop-gi}
For all smooth functions $f$ one has
$$
\forall x\in\dR^n,\,\,\forall t\geq 0,\,\,\frac{\partial }{\partial t}\Pt f(x)=\GI (\Pt f)(x)=\Pt (\GI f)(x),
$$
where for all smooth functions $f$, $\GI f=\Delta f-x\cdot\nabla f$.  

The linear operator $\GI$ is known as the {\it infinitesimal generator} of the Ornstein-Uhlenbeck semigroup.   
\end{eprop}
\begin{eproof}
If $f$ be a smooth function, then 
$$
\frac{\partial }{\partial t}\Pt f(x)=\int\PAR{-e^{-t}x+\frac{e^{-2t}}{\sqrt{1-e^{-2t}}}y}\cdot\nabla f\PAR{e^{-t}x+\sqrt{1-e^{-2t}}y}d\gamma(y). 
$$
By definition of the Ornstein-Uhlenbeck semigroup  one gets  
$$
-x{e^{-t}} \cdot\int \nabla f\PAR{e^{-t}x+\sqrt{1-e^{-2t}}y}d\gamma(y)=-x\cdot\nabla \Pt f(x)
$$ 
whereas the second term, after an integration by parts gives
$$
{\frac{e^{-2t}}{\sqrt{1-e^{-2t}}}}\int y\cdot\nabla f\PAR{e^{-t}x+\sqrt{1-e^{-2t}}y}d\gamma(y)=\Delta \Pt f(x), 
$$  
which finishes the proof. 

Using the same computation one can prove the commutation property between $\Pt$ and the generator $\GI$. 
\end{eproof}

More generally, if $\GI$ is an  infinitesimal generator associated to a linear semigroup $(\Pt)_{t\geq 0}$ (not necessary a Markov semigroup) then the commutation  $\GI\Pt=\Pt\GI$ holds. 

\begin{eprop}[Some properties of the O-U semigroup] 
\label{prop-er}The Ornstein-Uhlenbeck semigroup is $\gamma$-ergodic, that means for all $f\in\mathcal C_b(\dR^n)$, 
\begin{equation}
\label{eq-ergo}
\forall x\in\dR^n,\,\,\lim_{t\rightarrow\infty}\Pt f(x)=\int fd\gamma,
\end{equation}
in $L^2(d\gamma)$. 

The probability measure $\gamma$ is then the unique {\it invariant probability  measure}, for all smooth functions $f\in\mathcal C_b(\dR^n)$~: 
\begin{equation}
\label{eq-inv}
\int \Pt f d\gamma=\int f d\gamma, 
\end{equation}
or equivalently for all smooth functions $f$,
$$
\int \GI f d\gamma=0. $$ In fact we have the fundamental identity,
\begin{equation}
\label{eq-ipp} 
\int g\GI f d\gamma=\int f\GI g d\gamma=-\int\nabla f\cdot\nabla gd\gamma,
\end{equation}
for all smooth functions $f$ and $g$ on $\dR^n$. We say that   the Gaussian distribution is {\it reversible} with respect to the Ornstein-Uhlenbeck semigroup, $\GI$ is symmetric in $L^2(d\gamma)$.  
\end{eprop}
\begin{eproof}
Let us give the proof of~\eqref{eq-ipp}: 
$$
\begin{array}{rl}
\int f\GI gd\gamma 
&=\disp\int f\Delta gd\gamma-\int (fx\cdot\nabla g)d\gamma\\
&=-\disp\int \nabla\cdot(f\gamma)\cdot\nabla gdx-\int fx\cdot \nabla gd\gamma\\
&=-\disp\int \nabla f\cdot\nabla gd\gamma,
\end{array}
$$
where $\nabla\cdot f$ stands for the divergence of $f$.  

In fact~\eqref{eq-inv} is clear due to the fact if a semigroup is ergodic for some probability measure  then the measure is always invariant.  
\end{eproof}

As we have seen in the proof of Proposition~\ref{prop-gi}, the Ornstein-Uhlenbeck semigroup satisfies the equality for all $f$ and $x$:
\begin{equation}
\label{eq-comm}
\forall t\geq0,\,\,\nabla \Pt f(x)=e^{-t}\Pt \nabla f(x), 
\end{equation}
where $\Pt \nabla f= \PAR{\Pt \partial_i f}_{1\leq i\leq n}$ and for all norms $\NRM{\cdot}$ in $\dR^n$, one gets easily  
\begin{equation}
\label{eq-comm2}
\forall t\geq0,\,\,\NRM{\nabla \Pt f(x)}\leq e^{-t}\Pt \NRM{\nabla f}(x), 
\end{equation}
those equations are known as the commutation property of the gradient and the Ornstein-Uhlenbeck semigroup. Inequality~\eqref{eq-comm2}  is the key formula to get classical inequalities. 
\subsubsection{The Poincar\'e and logarithmic Sobolev inequalities}
\begin{ethm}
\label{thm-poin}
The following Poincar\'e inequality for the Gaussian measure holds, for all smooth functions $f$ on $\dR^n$, 
\begin{equation}
\label{eq-poin}
\var{\gamma}{f}:=\int f^2d\gamma-\PAR{\int fd\gamma}^2\leq\int\ABS{\nabla f}^2d\gamma.
\end{equation}

The term $\var{\gamma}{f}$ is {\it the variance of $f$ under  $\gamma$}. Moreover, the inequality is optimal and extremal functions are given by smooth functions satisfying $\nabla f=C$ for some constant   $C\in\dR^n$.  
\end{ethm}
\begin{eproof}
Let $f$ be a smooth function on $\dR^n$ then  $\PT{0} f=f$ and  $\PT{\infty} f=\int fd\gamma$ (see~\eqref{eq-ergo}), therefore the Ornstein-Uhlenbeck semigroup gives  a nice interpolation between $f$ and   $\int fd\gamma$. 
$$
\begin{array}{rl}
\disp\var{\gamma}{f}&=-\disp \int_0^{+\infty}\frac{d}{dt}\int \PAR{\Pt f}^2d\gamma dt \\
&=-2\disp\int_0^{+\infty}\int \GI \Pt f{\Pt f}d\gamma dt \\
&=2\disp\int_0^{+\infty}\int \ABS{\nabla\Pt f}^2d\gamma dt \\
&\leq \disp2\int_0^{+\infty}\int e^{-2t}(\Pt\ABS{\nabla f})^2d\gamma dt \\
&\leq2\disp\int_0^{+\infty}\int e^{-2t}\Pt\PAR{\ABS{\nabla f}^2}d\gamma dt \\
&=2\disp\int_0^{+\infty}\int e^{-2t}{\ABS{\nabla f}^2}d\gamma dt \\
&=\disp\int {\ABS{\nabla f}^2}d\gamma, \\
\end{array}
$$
where we use equality~\eqref{eq-comm2},  Cauchy-Schwarz inequality and the invariance property of the standard Gaussian distribution~\eqref{eq-inv}. 

On can check that in all stages of the proof, smooth functions satisfying $\nabla f=C$ are the unique function such that  the two inequalities become equalities.  
\end{eproof}

\begin{ethm}
\label{thm-logsob}
The following logarithmic Sobolev inequality  for the Gaussian measure holds,  for all smooth and non-negative functions $f$ on $\dR^n$, 
\begin{equation}
\label{eq-logsob}
\ent{\gamma}{f}:=\int f\log \frac{f}{\int fd\gamma}d\gamma\leq\frac{1}{2}\int\frac{\ABS{\nabla f}^2}{f}d\gamma. 
\end{equation}
The term $\ent{\gamma}{f}$ is known as {\it the entropy of $f$ under $\gamma$}. Moreover, the inequality~\eqref{eq-logsob} is optimal and extremal functions are given by 
$\nabla f=C f$ for some constant $C\in\dR^n$. 
\end{ethm}
\begin{eproof}
Let us mimic the proof of the Poincar\'e inequality, let $f$ be a smooth and non-negative function on $\dR^n$ then
$$
\begin{array}{rl}
\disp\ent{\gamma}{f}&=-\disp\int_0^{+\infty}\frac{d}{dt}\int {\Pt f}\log\Pt fd\gamma dt \\
&=-\disp\int_0^{+\infty}\int \GI \Pt f \log \Pt fd\gamma dt \\
&=\disp\int_0^{+\infty}\int {\nabla\Pt f}\cdot\nabla \log\Pt fd\gamma dt \\
&=\disp\int_0^{+\infty}\int \frac{\ABS{\nabla\Pt f}^2}{\Pt f}d\gamma dt, \\
&\leq \disp\int_0^{+\infty}\int e^{-2t}\frac{\PAR{\Pt \ABS{\nabla f}}^2}{\Pt f}d\gamma dt \\
\end{array}
$$
where we have used the same argument as for Poincar\'e inequality. Now Cauchy-Schwarz inequality  or the convexity of the map
$$
(x,y)\mapsto x^2/y
$$
for $x,y>0$, implies
$$
\frac{\PAR{\Pt \ABS{\nabla f}}^2}{\Pt f}\leq \Pt\PAR{\frac{\ABS{\nabla f}^2}{f}},
$$
 then one gets 
$$
\ent{\gamma}{f}\leq \int_0^{+\infty}\int e^{-2t}\Pt\PAR{\frac{\ABS{\nabla f}^2}{f}}d\gamma dt= \frac{1}{2}\int {\frac{\ABS{\nabla f}^2}{f}}d\gamma.
$$
One obtains extremal functions in the same way than for Poincar\'e inequality. 
\end{eproof}

The logarithmic Sobolev inequality is often noted for  $f^2$ instead of $f$, which gives for all smooth functions $f$,
$$
\ent{\gamma}{f^2}\leq 2\int{\ABS{\nabla f}^2}d\gamma.
$$
{\it At the light of the Theorems~\ref{thm-poin} and \ref{thm-logsob}, we say that  the standard  Gaussian satisfies  a Poincar\'e and a  logarithmic Sobolev inequality.  }

More generally a logarithmic Sobolev inequality always implies a Poincar\'e inequality by a Taylor expansion (see Chapter 1 of \cite{logsob}).

In proposition~\ref{prop-er}, we proved that the Ornstein-Uhlenbeck semigroup is ergodic with respect to the Gaussian distribution. In fact  one of the main application of the Poincar\'e and the logarithmic Sobolev inequalities is to give an estimate of the speed of convergence in two different spaces. 
\begin{ethm}
\label{thm-speed}
The Poincar\'e inequality~\eqref{eq-poin} is equivalent to the following inequality 
\begin{equation}
%\label{eq-l2}
\var{\gamma}{\Pt f}\leq e^{-2t}\var{\gamma}{f},
\end{equation}
for all smooth functions $f$. 

And in the same way, the logarithmic Sobolev inequality~\eqref{eq-logsob} is equivalent to 
\begin{equation}
%\label{eq-llogl}
\ent{\gamma}{\Pt f}\leq e^{-2t}\ent{\gamma}{f},
\end{equation} 
for all non-negative and smooth  functions $f$. 
\end{ethm}

\begin{eproof}
For the first assertion, an elementary computation gives that 
$$
\frac{d}{dt}\var{\gamma}{\Pt f}=-2\int\ABS{\nabla \Pt f}^2d\gamma,
$$  
then the Poincar\'e inequality and  Gr\"onwall lemma implies~\eqref{eq-l2}. Conversely, the derivation at time $t=0$ of~\eqref{eq-l2} implies the Poincar\'e inequality.   

For the second assertion, we use the same method and the derivation of the entropy, 
\begin{equation}
\label{eq-deri}
\frac{d}{dt}\ent{\gamma}{\Pt f}=-\int\frac{\ABS{\nabla \Pt f}^2}{\Pt f}d\gamma. 
\end{equation}
\end{eproof}

One of the main difference between the two inequalities is that the initial condition is in $L^2(d\gamma)$ for the Poincar\'e inequality whereas the initial condition is in $L\log L(d\gamma)$ for the logarithmic Sobolev inequality.  
%%%%%%%%%%%%
%%%%%%%%%%%%
\section{Poincar\'e and logarithmic Sobolev inequalities under curvature criterium}
\label{curvature}
%%%%%%%%%%%%
%%%%%%%%%%%%

The main idea of this section is to obtain criteria for a probability measure $\mu$ such that the two inequalities~\eqref{eq-poin} and~\eqref{eq-logsob} hold for the measure $\mu$. We will study a particular case of  the curvature-dimension criterium (or $\Gamma_2$-criterium) introduced by D. Bakry and M. Emery in \cite{bakryemery}. This criterium gives conditions on an infinitesimal generator $\GI$ such that all the  computations done for the Ornstein-Uhlenbeck semigroup could be applied to $\GI$. 

Let a function $\psi\in \mathcal C^2(\dR^n)$, and define the infinitesimal generator: 
\begin{equation}
\label{eq-gi}
\GI f=\Delta f-\nabla \psi\cdot\nabla f,
\end{equation}
for all smooth functions $f$. 

Assume that $\int e^{-\psi}dx<+\infty$ and define the probability measure  $d\mu_\psi(x)=\frac{ e^{-\psi}dx}{Z_\psi}dx$, where $Z_\psi=\int e^{-\psi}dx$. It is easy to see that the operator $\GI$ satisfies for all smooth functions $f$ and $g$ on $\dR^n$,  
\begin{equation}
\label{eq-ipp2}
\int f\GI gd\mu_\psi=\int g\GI fd\mu_\psi =-\int \nabla f\cdot\nabla gd\mu_\psi,
\end{equation}
and $
\int \GI fd\mu_\psi=0$. We recover the same property as for the Ornstein-Uhlenbeck semigroup, see~\eqref{eq-ipp}. 
The generator $\GI$ is symmetric in $L^2(d\mu_\psi)$ and the probability measure $\mu_\psi$ is also invariant with respect to $\GI$. 

Let define the {\it Carr\'e du champ}, for all smooth functions $f$, 
\begin{equation}
\label{eq-cc}
\Gamma(f,f)=\frac{1}{2}\PAR{\GI (f^2)-2f\GI f},
\end{equation}
we note usually $\Gamma(f)$ instead of $\Gamma(f,f)$. 
The carr\'e du champ is a quadratic form and  the bilinear form associated is given by 
$$
\Gamma(f,g)=\frac{1}{2}\PAR{\GI (fg)-f\GI g-g\GI f}.
$$
If we iterate the process one obtains   the $\Gamma_2$-operator, for all smooth functions $f$, 
\begin{equation}
\label{eq-gamma2}
\Gamma_2(f,f)=\frac{1}{2}\PAR{\GI (\Gamma(f))-2\Gamma(f,\GI f)}.
\end{equation}

We assume in this section that there exits a set  of function $\A$, dense in $L^2(d\mu)$, such that all computations can be done in this class of function. In the previous section, the set $\A$ was $\mathcal C_c^\infty(\dR^n)$ and one of the main problem is to describe this class of functions. It can be done under the $\Gamma_2$-criterium $CD(\rho,+\infty)$ (see the definition  below),  we refer to~\cite{logsob,bakrytata} and references therein to get  more  informations.

\begin{edefi}
We say that the linear operator $\GI$, satisfies the $\Gamma_2$-criterium $CD(\rho,+\infty)$ with some $\rho\in\dR$, if for all functions $f\in\A$ 
\begin{equation}
\label{def-gamma22}
\Gamma_2\PAR{f}\geq\rho \Gamma(f). 
\end{equation}
\end{edefi}

\begin{erem}
Since for all smooth functions $f$, $\GI f=\Delta f-\nabla \psi\cdot\nabla f$, a straight forward computation gives, 
$$
\Gamma(f)=\ABS{\nabla f}^2,
$$ 
and 
$$
\Gamma_2(f)=\NRM{\rm{Hess}(f)}_{H.S.}^2+<\nabla f,\rm{Hess}(\psi)\nabla f>,
$$
where the Hilbert-Schmidt norm is given by   $\NRM{\rm{Hess}(f)}_{H.S.}^2=\sum_{i,j}\PAR{\frac{\partial^2}{\partial x_i\partial x_j}f }^2$.

Then the linear operator $\GI$ defined in~\eqref{eq-gi}  satisfies the $\Gamma_2$-criterium $CD(\rho,+\infty)$ with some $\rho\in\dR$ if for all $x\in\dR^n$
\begin{equation}
\label{def-gamma2}
\rm{Hess}({\psi})(x)\geq\rho Id,
\end{equation}
in the sense of the symmetric matrix, i.e. for all $Y\in\dR^n$,
$$
<Y,\rm{Hess}({\psi})(x)Y>\,\geq\rho \ABS{Y}^2,
$$
where $<\cdot,\cdot>$ is the Euclidean scalar product.  
\end{erem}

\begin{ethm}
\label{thm-gamma2}
Let $\psi\in\mathcal C^2(\dR^n)$ and assume that there exists $\rho>0$ such that the linear operator~\eqref{eq-gi} satisfies a $\Gamma_2$-criterium $CD(\rho,+\infty)$, then the probability measure $\mu_\psi$ satisfies a Poincar\'e inequality   
\begin{equation}
\label{eq-poin2}
\var{\mu_\psi}{f}\leq\frac{1}{\rho}\int\ABS{\nabla f}^2d\mu_\psi,
\end{equation}
for all $f\in\A$ and a logarithmic Sobolev inequality
\begin{equation}
\label{eq-logsob2}
\ent{\gamma}{f}\leq\frac{1}{2\rho}\int\frac{\ABS{\nabla f}^2}{f}d\mu_\psi,
\end{equation}
for all smooth and non-negative functions $f\in\A$. 
\end{ethm} 
\begin{elem}
Let $(\Pt)_{t\geq0}$ be the Markov semigroup associated to the infinitesimal generator $\GI$. Assume that $\rho>0$ then  $(\Pt)_{t\geq0}$ is $\mu_\psi$-ergodic which means for all functions $f\in\A$ 
$$
\lim_{t\rightarrow+\infty}\Pt f (x)=\int fd\mu_\psi,
$$
in $f\in L^2(d\mu_\psi)$ and $\mu_\psi$ almost surely. 
\end{elem}
\begin{elem}
\label{lem-2}
Let $\phi$ be a $\mathcal C^2$ function, then for all functions $f\in\A$, 
\begin{equation}
\label{eq-dif1}
\GI \phi(f)=\phi'(f)\GI f+\phi''(f)\Gamma(f)\,\,\rm{and}\,\,\Gamma(\log f)=\frac{1}{f^2}\Gamma(f),
\end{equation}
moreover one has 
\begin{equation}
\label{eq-dif2}
\Gamma_2(\log f)=\frac{1}{f^2}\Gamma_2(f)-\frac{1}{f^3}\Gamma(f,\Gamma(f))+\frac{1}{f^4}\PAR{\Gamma(f)}^2
\end{equation}
\end{elem}

{\noindent {\emph{\textbf{Proof of the Theorem~\ref{thm-gamma2} }}}\\$\lhd$~}
First we prove the first inequality~\eqref{eq-poin2}. As for the Ornsten-Uhlenbeck semigroup, one gets if $(\Pt)_{t\geq0}$ is the Markov semigroup associated to the infinitesimal generator $\GI$, for all functions $f\in\A$,
$$
\begin{array}{rl}
\var{\mu_\psi}{f}&=-\disp\int_0^{+\infty}\frac{d}{dt}\int \PAR{\Pt f}^2d\mu_\psi dt \\
&=-2\disp\int_0^{+\infty}\int \GI \Pt f{\Pt f}d\mu_\psi dt 
\end{array}
$$
Since $\mu_\psi$ is invariant,
$$
\int2\Pt f\GI\Pt fd\mu_\psi=\int\PAR{2\Pt f\GI\Pt f-\GI(\Pt f)^2 }d\mu_\psi=-2\int\Gamma(\Pt f)d\mu_\psi, 
$$
which gives  
\begin{equation}
\label{eq-poinc}
\var{\mu_\psi}{f}=\int_0^{+\infty}2\int\Gamma(\Pt f)d\mu_\psi dt.
\end{equation}
 Let now consider for all $t>0$, 
$$
\Phi(t)=2\int\Gamma(\Pt f)d\mu_\psi,
$$ 
The time derivative of $\Phi$ is equal to 
\begin{multline*}
\Phi'(t)=4\int\Gamma(\Pt f,\GI\Pt f)d\mu_\psi=\\2\int \PAR{ 2\Gamma(\Pt f,\GI\Pt f)-\GI(\Gamma(\Pt f)) }d\mu_\psi=-4\int\Gamma_2(\Pt f)d\mu_\psi.
\end{multline*}
The $\Gamma_2$-criterium implies that  $\Phi'(t)\leq -2\rho\Phi(t)$ which gives $\Phi(t)\leq e^{-{t2}{\rho}}\Phi(0)$. The last inequality with~\eqref{eq-poinc} implies 
$$
\var{\mu_\psi}{f}\leq\int_0^{+\infty}e^{-{t2}{\rho}}dt\int 2\Gamma(f)d\mu_\psi=\frac{1}{\rho}\int \Gamma(f)d\mu_\psi dt.
$$
Let now prove the logarithmic Sobolev inequality for the measure $\mu_\psi$. Let $f$ be a non-negative and smooth function on $\dR^n$,  
$$
\begin{array}{rl}
\ent{\mu_\psi}{f}&=-\disp\int_0^{+\infty}\frac{d}{dt}\int {\Pt f}\log \Pt fd \mu_\psi dt \\
&=-\disp\int_0^{+\infty}\int \GI \Pt f\log {\Pt f}d\mu_\psi dt 
\end{array}
$$
Since $\GI$ is symmetric and by lemma~\ref{lem-2} one gets 
$$
\int \GI \Pt f\log{\Pt f} d\mu_\psi =\int \Pt f \GI\log{\Pt f} d\mu_\psi= -\int \frac{\Gamma(\Pt f)}{\Pt f}d\mu_\psi=-\int {\Gamma(\log \Pt f)}{\Pt f}d\mu_\psi,
$$
which gives  
\begin{equation}
\label{eq-logo}
\ent{\mu_\psi}{f}=\int_0^{+\infty}\int {\Gamma(\log\Pt f)}{\Pt f}d\mu_\psi dt.
\end{equation}
As for Poincar\'e inequality, let consider for all $t>0$, 
$$
\Phi(t)=\int\frac{\Gamma(\Pt f)}{\Pt f}d\mu_\psi
$$ 
where $\Pt f=g$. 
The time derivative of $\Phi$ is equal to 
$$
\Phi'(t)=\int \PAR{2\frac{\Gamma(\GI g,g)}{g}-\frac{\GI g \Gamma(g)}{g^2}}\mu_\psi=\int \PAR{2\frac{\Gamma(\GI g,g)}{g}-\frac{\GI g \Gamma(g)}{g^2}-\GI\PAR{\frac{\Gamma(g)}{g}}}\mu_\psi.
$$
Since $$
\GI\PAR{\frac{\Gamma(g)}{g}}=2\Gamma\PAR{\Gamma(g),\frac{1}{g}}+\frac{1}{g}\GI\Gamma(g)+\GI \PAR{\frac{1}{g}}\Gamma(g),
$$
by Lemma~\ref{lem-2} one has 
$$
\Phi'(t)=-2\int\Gamma_2(\log \Pt f)\Pt fd\mu_\psi. 
$$
The $\Gamma_2$-criterium implies that  $\Phi'(t)\leq -2\rho\Phi(t)$ which gives $\Phi(t)\leq e^{-2{\rho t}}\Phi(0)$. This inequality  with~\eqref{eq-logo} implies that
$$
\ent{\mu_\psi}{f}\leq\int_0^{+\infty}e^{-{}{2\rho t}}dt\int \Gamma (\log f)fd\mu_\psi=\frac{1}{2\rho}\int \Gamma(\log f)fd\mu_\psi=\frac{1}{2\rho}\int \frac{\ABS{\nabla f}}{f}d\mu_\psi.
$$
{~$\vartriangleright$\\}

The meaning of this result is : if $\mu_\psi$ is more $\log$-concave than the Gaussian distribution then $\mu_\psi$ satisfies both inequalities. 
\begin{erem}
The $\Gamma_2$-criterium is in fact a more general criterium. The definition of a  diffusion semigroup could be a Markov semigroup such that for all smooth functions $\phi$, the  equations~\eqref{eq-dif1} and~\eqref{eq-dif2} hold for the generator associated to the semigroup.

 In fact on $\dR^n$ (or on a manifold on a local chart)  that means that the infinitesimal generator $\GI$ of the Markov semigroup is given by,
$$
\forall x\in\dR^n,\,\,\GI f(x) = \sum_{i,j}D_{i,j}(x)\partial_{i,j} f(x)-\sum_i a_i(x)\partial_i f (x),
$$
where $D(x)=(D_{i,j}(x))_{i,j}$ is a symmetric and non-negative matrix and $a(x)=(a_i(x))_{i}$ is a vector. 

Then the conditions $\Gamma_2(f)\geq\rho \Gamma(f)$ for some $\rho>0$ implies that there exists an invariant measure $\mu$ of the semigroup and $\mu$ satisfies the   Poincar\'e and a logarithmic Sobolev inequality with the same constant as before. One of the difficulties of this general case is to find tractable conditions on functions $D$ and $a$ such that the $\Gamma_2$-criterium holds. Some others examples can be found in~\cite{bolley-gentil}. 

Let us also note that the $\Gamma_2$-criterium $CD(\rho, \infty)$ is a particular case of the $CD(\rho,n)$ criterium where $n\in\dN^*$ :
$$
\Gamma_2(f)\geq\rho\Gamma(f)+\frac{1}{n}\PAR{\GI f}^2,
$$  
for all smooth functions $f$. For example, the Ornstein-Uhlenbeck semigroup satisfies the $CD(1,\infty)$ criterium and  the heat equation $\GI=\Delta$ satisfies the $CD(0,n)$.   On can observe that the Ornstein-Uhlenbeck semigroup does not satisfies a $CD(r,m)$ criterium for any  $r,m>0$. 
\end{erem}

\begin{ethm}
\label{thm-speed-gene} As for the Ornstein-Uhlenbeck semigroup, 
the Poincar\'e inequality~\eqref{eq-poin2} is equivalent to the following inequality 
\begin{equation}
\label{eq-l2}
\var{\mu_\psi}{\Pt f}\leq e^{-\frac{2}{\rho}t}\var{\mu_\psi}{f},
\end{equation}
for all  functions $f\in \A$. 

And in the same way, the logarithmic Sobolev inequality~\eqref{eq-logsob2} is equivalent to 
\begin{equation}
\label{eq-llogl}
\ent{\mu_\psi}{\Pt f}\leq e^{-2t}\ent{\mu_\psi}{f},
\end{equation} 
for all non-negative functions $f\in\A$. 
\end{ethm}

The logarithmic Sobolev inequality has two main applications. The first one the asymptotic behaviour in term of entropy, this is the result of Theorem~\ref{thm-speed-gene}. The second application is about concentration inequality, a probability measure $\mu$ satisfying a logarithmic Sobolev inequality has the same tail as the Gaussian distribution.

 This properties is also a consequence of  the Talagrand inequality described in the next section. 
%%%%%%%%%%%%
%%%%%%%%%%%%
\section{The logarithmic Sobolev  and transport inequalities by transportation method}
\label{sec-tr}
%%%%%%%%%%%%
%%%%%%%%%%%%

%\subsection{Mass transportation method}
%\label{int}
We will see   how Brenier's Theorem  can be used in this context to give a new proof of the logarithmic Sobolev inequality, the method is called mass transportation method.   

We will illustrate this method for the Gaussian measure but it could be generalized for a large class of measures, this will be discussed later.  The method come from~\cite{ov,c} and has been generalized for many Euclidean inequalities as Sobolev and Gagliardo-Nirenberg inequalities, see~\cite{agk,cnv,nazaret}.  

\medskip

The Wasserstein distance between two probability measures $\mu$ and $\nu$ is defined by  
\begin{eqnarray}
\label{eqdistance}
 W_2(\mu,\nu)=\left(\inf \int \ABS{x-y}^2d\pi(x,y) \right)^{1/2}.
\end{eqnarray}
where the infimum is running over all probability measures $\pi$ on $\dR^n\times\dR^n$ with respective marginals $\mu$ and $\nu$~: for all  bounded functions $g$ and $h$,  
$$
\int (g(x)+h(y))d\pi(x,y)=\int gd\mu+\int hd\nu.
$$
Such probability is called a coupling of $(\mu,\nu)$. 

\medskip

Brenier's theorem says that  that there exits an optimal deterministic  coupling of $(\mu,\nu)$ : there exists  a convex  map $\Phi$ satisfying
\begin{equation*}
\int h(\nabla \Phi)d\nu=\int hd\mu, 
\end{equation*}
for all bounded  functions $h$. Moreover 
$$
W_2^2(d\nu,d\mu)=\int\ABS{\nabla \theta}^2d\nu,
$$
where  $\theta(x) = \Phi(x)-\frac{1}{2}\ABS{x}^2$.  This result has been proved by Brenier, $\nabla\Phi$ is called the Brenier map between $\nu$ and $\mu$, see~\cite{villani}. 

\medskip

We apply this result in the Gaussian case. Let $f$ be a smooth and positive  function  such that $\int fd\gamma = 1$, Brenier's theorem implies that there exists  a convex  map $\Phi$ satisfying
\begin{equation}
\label{eq-chang}
\int h(\nabla \Phi)fd\gamma=\int hd\gamma, 
\end{equation}
for all bounded and measurable functions $h$. Moreover 
$$
W_2^2(fd\gamma,d\gamma)=\int\ABS{\nabla \theta}^2fd\gamma,
$$
where  $\theta(x) = \Phi(x)-\frac{1}{2}\ABS{x}^2$.

\medskip

If now $\Phi$ is a $\mathcal C^2(\dR^n)$ function, then coming from~\eqref{eq-chang}, the Monge-Amp\`ere equation holds~:
$fd\gamma$-a.e. 
\begin{equation}
\label{eq-ma}
f (x)e^{-{|x|^2}/{2}} = \rm{det}\PAR{Id + \rm{Hess}(\theta)} e^{ - |x+\nabla\theta(x)|^2/2}.
\end{equation}
After taking the logarithm, we get
$$
\begin{array}{rl}
\log f (x) &=  - \disp\frac{1}{2} \ABS{x +\nabla\theta(x)}^2 + \frac{1}{2}|x|^2 + \log \det(Id + \rm{Hess}(\theta))\\
&=  -\disp x\cdot\nabla\theta(x)  -  \frac{1}{2} \ABS{\nabla\theta(x)}^2 + \log \det(Id + \rm{Hess}( \theta))\\
&\leq  - \disp x\cdot\nabla\theta(x)  - \frac{ 1}{2} \ABS{\nabla\theta(x)}^2 + \Delta\theta(x),
\end{array}
$$
where we used inequality $\log(1 + t)\leq  t$ whenever $1 + t>   0$. We integrate with respect to
$fd\gamma$ :
$$
\ent{\gamma}{f}\leq \int f\PAR{\Delta \theta-x\cdot\nabla\theta}d\gamma-\int \frac{1}{2} \ABS{\nabla\theta(x)}^2fd\gamma.
$$
The integration by parts implies 
$$
\begin{array}{rl}
\ent{\gamma}{f}&\leq -\disp\int \nabla \theta\cdot\nabla fd\gamma-\int \frac{1}{2} \ABS{\nabla\theta(x)}^2fd\gamma\\
& \leq-\disp\frac{1}{2}\int \ABS{\sqrt{f}\nabla \theta+\frac{\nabla f}{\sqrt{f}}}^2d\gamma+\frac{1}{2}\int \frac{\ABS{\nabla f}^2}{f}d\gamma\\
& \leq\disp\frac{1}{2}\int \frac{\ABS{\nabla f}^2}{f}d\gamma,\\
\end{array}
$$
which is the optimal logarithmic Sobolev inequality~\eqref{eq-logsob}. 

\medskip

Hence we have proved, using Brenier's map,  the logarithmic Sobolev inequality for the Gaussian measure with the optimal constant. As we can see in the proof, one has assumed that $\Phi$ is a $\mathcal C^2$ function. It can be obtained using Caffarelli's regularity theory :  it needs another assumptions, $f$ has to be  smooth with a  compact and convex support. We skip it for simplicity of the description of the method, many informations can be bound in~\cite{villani}

\bigskip

Let us see what can be done if now $\nabla \Phi$ be the Brenier map between $d\gamma$ and $fd\gamma$ instead $fd\gamma$ and $d\gamma$ :  that is for all bounded and measurable functions $h$: 
$$
\int hfd\gamma=\int h(\nabla \Phi)d\gamma,
$$
and if  $x+\nabla\theta(x)=\nabla \Phi$ then 
$$
W_2^2(fd\gamma,d\gamma)=\int\ABS{\nabla \theta}^2d\gamma.
$$
In that case the Monge-Amp\`ere equation gives 
\begin{equation}
\label{eq-ma2}
\rm{det}\PAR{Id + \rm{Hess}(\theta)}f (x+\nabla\theta(x))e^{-{|x+\nabla\theta(x)|^2}/{2}} =  e^{ - |x|^2/2}.
\end{equation}
Which implies
$$
\begin{array}{rl}
\log f (x+\nabla\theta(x)) &=  \disp\frac{1}{2} \ABS{x +\nabla\theta(x)}^2 - \frac{1}{2}|x|^2 - \log \det(Id + \rm{Hess}(\theta))\\
&=   \disp x\cdot\nabla\theta(x)  +  \frac{1}{2} \ABS{\nabla\theta(x)}^2 - \log \det(Id + \rm{Hess}( \theta))\\
&\geq  \disp x\cdot\nabla\theta(x) +  \frac{ 1}{2} \ABS{\nabla\theta(x)}^2 - \Delta\theta(x)\\
&\disp =-\GI\theta +  \frac{ 1}{2} \ABS{\nabla\theta(x)}^2,
\end{array}
$$
where $\GI$ is the Ornstein-Uhlenbeck generator. Then 
$$
\begin{array}{rl}
\ent{\gamma}{f} &= \disp\int f\log fd\gamma\\
&=\disp\int \log f(\nabla \Phi)d\gamma\\
&\geq \disp\int -L\theta d\gamma+  \int \frac{ 1}{2} \ABS{\nabla\theta(x)}^2d\gamma \\
&= \int \disp\frac{ 1}{2} \ABS{\nabla\theta(x)}^2d\gamma =\disp\frac{ 1}{2}W_2^2(fd\gamma, d\gamma)
\end{array}
$$
We have  proved that for all functions $f$ such that $fd\gamma$ is a probability measure, one has 
\begin{equation}
\label{eq-t2}
W_2(fd\gamma, d\gamma)\leq \sqrt{2\ent{\gamma}{f}}.
\end{equation}
This inequality, called {\it Talagrand inequality for the Gaussian distribution   (or $\mathcal T_2$ inequality)}, has been proved by Talagrand in~\cite{t1}. 

As for Poincar\'e and logarithmic Sobolev inequalities,  we says that a probability measure $\mu$ satisfies a Talagrand inequality  if there exists $C\geq 0$ such that,
\begin{equation}
\label{eq-t22}
T_2(fd\mu, d\mu)\leq \sqrt{C\ent{\mu}{f}},
\end{equation}
for all functions $f$ such that $fd\mu$ is a probability measure,

\subsection{Remarks and extensions}
This method can also be used is the  context of the section~\ref{curvature}.  Assume that $\psi$ is uniformly convex and satisfying 
$$
\rm{Hess}(\psi)\geq \rho I, 
$$
with some $\rho>0$. The mass transportation method  implies that the measure 
$$
d\mu_\psi(x)=\frac{ e^{-\psi}dx}{Z_\psi}dx
$$
 satisfies the logarithmic Sobolev inequality~\eqref{eq-logsob2} with the constant $1/(2\rho)$. This is an alternative proof of Theorem~\ref{thm-gamma2}. Actually this method is not useful to obtain directly a Poincar\'e inequality. 

\medskip

Of course, as for Ornstein-Uhlenbeck semigroup,  the mass transportation method gives also a Talagrand inequality~\eqref{eq-t22} for the measure $\mu_\psi$ :  
$$
T_2(fd\mu_\psi, d\mu_\psi)\leq \sqrt{\frac{1}{\rho}\ent{\mu_\psi}{f}},
$$
for all probability measure $fd\mu_\psi$. 

\medskip

In fact the general result holds, 
\begin{ethm}[Otto-Villani]
Let $\mu$ be a probability measure on $\dR^n$ satisfying a logarithmic Sobolev inequality
 $$
\ent{\mu}{f^2}\leq C\int{\ABS{\nabla f}^2}d\mu,
$$
for all smooth functions $f$ and for some constant $C\geq0$. 

Then $\mu$ satisfies a Talagrand inequality 
$$
T_2(fd\mu, d\mu)\leq \sqrt{2C\ent{\mu}{f}},
$$
for all probability measure $fd\mu$.  
\end{ethm}
The original proof comes from  \cite{ov} and an easier one,  using Hamilton-Jacobi equation, has been given in~\cite{bgl}. These two inequalities are quite similar but it has been proved in~\cite{cagu,goz} that they are not equivalent.

\newcommand{\etalchar}[1]{$^{#1}$}

\bigskip
\noindent

\medskip\noindent

\noindent
Ceremade, UMR CNRS 7534 \\
Universit\'e Paris-Dauphine\\
Place du Mar\'echal De Lattre De Tassigny\\
75116 PARIS - France\\
gentil@ceremade.dauphine.fr\\

\end{document}